\documentclass[graybox]{svmult}
\usepackage{amsmath,amsfonts,amssymb,mathrsfs} 
\usepackage{graphicx}
\newcommand*{\abs}[1]{\left\vert #1\right\vert}
\newcommand*{\dif}{\mathrm{d}}

\newcommand*{\e}{\mathrm{e}}
\newcommand*{\mi}{\mathrm{i}}

\newcommand*{\cnj}[1]{\overline{#1}}
\newcommand*{\hardy}[1]{\mathsf{H}^{#1}}
\newcommand*{\hardyT}[1]{\mathsf{H}^{#1}({\mathbb T})}
\newcommand*{\hardyR}[1]{\mathsf{H}^{#1}({\mathbb R})}
\newcommand*{\uhp}{{\mathbb H}}
\newcommand*{\blp}{{\mathsf B}}
\newcommand*{\mat}[1]{{\mathsf #1}}
\newcommand*{\radon}[2]{{\mathsf R}_{#1}{#2}}
\newcommand*{\ind}{{\large\bf 1}}
\newcommand*{\Dif}{{\mathsf D}}

\begin{document}

\title*{ On Complex Analytic tools, and the Holomorphic Rotation methods}
\author{Ronald R.~Coifman, Jacques Peyri\`ere, and Guido Weiss}
\institute{Ronald R.~Coifman \at Department of Mathematics, Program in
  Applied Mathematics, Yale University, New Haven, CT 06510, USA,
  \email{coifman-ronald@yale.edu} \and Jacques Peyri\`ere \at Institut
  de Math\'ematiques d'Orsay, CNRS, Universit\'e Paris-Saclay, 91405
  Orsay, France, \email{jacques.peyriere@universite-paris-saclay.fr}}



\maketitle

\section{introduction}             
  This paper in honor of Guido Weiss was written posthumously, jointly with him, as we had, all of his initial notes and ideas related to the program described below.
  
   Our task, here, is to recount ideas, explorations, and visions that Guido  his collaborators and students, developed over the last 60 years.  To point out the connection of ideas between the original views of the interplay between complex and real analysis  as envisioned by Zygmund and his students Calder\'on, Guido Weiss, Eli Stein,and many others,  70 years ago, and the current  approaches introducing  nonlinear multi layered analysis for the organization and processing of complicated oscillatory functions.   
    
     It was Zygmund’s view that harmonic analysis provides the infrastructure linking most areas of analysis, from complex analysis to partial differential equations, to probability, number theory, and geometry.
    
     In particular he pushed forward the idea that the remarkable tools of complex analysis, which include; contour integration, conformal mappings, factorization, tools which were used to provide miraculous proofs in real analysis, should be deciphered and converted to real variable tools.
 Together with Calder\'on, they bucked the trend for abstraction, prevalent at the time, and  formed a school pushing forward this interplay between real and complex analysis.
   A principal bridge was provided by real variable methods, multiscale analysis, Littlewood Paley theory, and related  Calderon representation formulas. Our aim, here, is to elaborate on the "magic" of complex analysis and indicate potential applications in Higher dimensions.
  An old idea of Calder\'on and Zygmund, the so called "rotation method", enabled the reductions of the study of $L^p$ estimates for multi dimensional singular integrals to a superposition,over all directions,of Hilbert transforms. Thereby allowing the use of one complex variable methods. A related idea was the invention of systems of Harmonic functions satisfying generalised Cauchy Riemann equations, such as the Riesz systems, exploiting their special properties. \cite{CW} 
  
  Our goal is to extend these ideas to enable remarkable nonlinear complex analytic tools for the adapted analysis of functions in one variable, to apply in higher dimensions.

  Guido has been pushing the idea that factorization theorems like Blaschke products are a key to  a variety of nonlinear analytic methods ~\cite{CW1}. Our goal here is to demonstrate this point, deriving amazing approximation theorems, in one variable, and opening doors to higher dimensional applications. Application in which each harmonic function is the average of special holomorphic functions in  planes and constant in orthogonal directions.

  We start by describing  recent developments in nonlinear complex analysis, exploiting the tools of factorization and composition. In particular we will sketch methods extending conventional Fourier analysis, exploiting both phase and amplitudes of holomorphic functions. 
  The "miracles of nonlinear complex analysis", such as factorization and composition of functions lead to new versions of holomorphic wavelets, and  relate them to multiscale dynamical systems.
  
  Our story  interlaces the role of the phase of signals with their analytic/geometric properties.  The Blaschke factors are a key ingredient, in building analytic tools, starting with the Malmquist-Takenaka orthonormal bases of the Hardy space $\hardyT{2}$, continuing with "best" adapted bases obtained through phase unwinding, and describing relations to composition of Blaschke products and their dynamics (on the disc and upper half plane).
  Specifically we construct multiscale orthonormal holomorphic wavelet bases,  generalized scaled holomorphic orthogonal bases, to dynamical systems, obtained by composing Blaschke products.  
  
   We also, remark, that the phase of a Blaschke product is a one layer neural net with ($\arctan$ as an activation sigmoid) and that the composition is a "Deep Neural Net" whose "depth" is the number of compositions. Our results provide a wealth of related libraries of orthogonal bases.
  
  We  sketch these ideas in various "vignette" subsections and
 refer for more details on analytic methods~\cite{CP}, related to the Blaschke based nonlinear phase 
unwinding decompositions~\cite{coifman,CSW,nahon}. We also consider
orthogonal decompositions of invariant subspaces of Hardy spaces. In
particular we constructed a multiscale decomposition, described below, of the Hardy space of the upper half-plane. 

Such a decomposition can
be carried in the unit disk by conformal mapping. A somewhat different
multiscale decomposition of the space $\hardyT{2}$ has been
constructed by using Malmquist-Takenaka bases associated with Blaschke
products whose zeroes are
$\displaystyle (1-2^{-n})\e^{2\mi\pi j/2^n}$ where $n\ge 1$ and
$0\le j< 2^n$ \cite{feichtinger}. Here we provide a variety  of multiscale decompositions by considering iterations of Blaschke products.

In the next chapter we will show how with help of an extended Radon transform we can introduce  a method of rotations to enable us to lift the one dimensional tools to higher dimensions. In particular the various orthogonal bases of holomorphic functions in one dimension, give rise to orthogonal bases of Harmonic functions in the higher dimensional upper half space.

\section{Preliminaries and notation}

For $p\ge 1$, $\hardyT{p}$ stands for the space of analytic functions~$f$
on the unit disk~${\mathbb D}$ such that
\begin{equation*}
\sup_{0< r<1} \int_{0}^{2\pi} |f(r\e^{\mi\theta})|^p\frac{\dif
  \theta}{2\pi} < +\infty.
\end{equation*}
Such functions have boundary values almost everywhere, and the Hardy
space $\hardyT{p}$ can be identified with the set of $L^p$ functions
on the torus~${\mathbb T}=\partial{\mathbb D}$ whose Fourier
coefficients of negative order vanish. We will alternate between analysis on the disk, and
the  parallel theory for analytic functions on the upper half
plane $\uhp = \{x+\mi y \ :\ y>0\}$. The space of analytic
functions~$f$ on $\uhp$ such that
$$\sup_{y>0} \|f(\cdot+\mi y)\|_{L^p({\mathbb R})} < +\infty$$
is denoted by $\hardyR{p}$. These functions have boundary values in
$L^p({\mathbb R})$ when $p\ge 1$. The space $\hardyR{p}$ is identified
to the space of $L^p$ functions whose Fourier transform vanishes on
the negative half line~$(-\infty,0)$.

\section{ Analysis on The upper half plane}

We present some known results~\cite{CP}, without proof.
 In this section one simply writes $\hardy{2}$ instead of
$\hardyR{2}$.

\subsection*{Malmquist-Takenaka bases}

Let
$(a_j)_{1\le j}$ be a sequence (finite or not)) of complex numbers
with positive imaginary parts and such that
\begin{equation}\label{upper} 
\displaystyle \sum_{j\ge 0} \frac{\Im a_j}{1+|a_j|^2} < +\infty.
\end{equation}
The corresponding Blaschke product is
$$\blp(x) = \prod_{j\ge 0}
\frac{\abs{1+a_j^2}}{1+a_j^2}\,\frac{x-a_j}{x-\cnj{a}_j},$$
where, $0/0$, which appears if $a_j=\mi$, should be understood as~1.
The factors $\displaystyle \frac{\abs{1+a_j^2}}{1+a_j^2}$ insure the
convergence of this product when there are infinitely many
zeroes. But, in some situations, it is more convenient to use other
convergence factors as we shall see below.

Whether the series~\eqref{upper} is convergent or not, one defines
(for $n\ge 0$) the functions
\begin{equation*}
\phi_n(x) = \frac{1}{\sqrt{\pi}}\left( \prod_{0\le j< n}
\frac{x-a_j}{x-\cnj{a}_j}\right)\, \frac{1}{x-\cnj{a}_n}.
\end{equation*}
Then these functions form an orthonormal system in $\hardy{2}$. If the
series~\eqref{upper} diverges, it is a  Malmquist-Takenaka orthonormal basis of $\hardy{2}$, otherwise
it is a basis of the orthogonal complement of $\blp\,\hardy{2}$ in
$\hardy{2}$. 

We  remark that roughly a hundred years ago these bases were constructed~\cite{takenaka,malmquist} through a Gram Schmidt orthogonalization of the  list of rational functions with poles in the lower half plane .  

Observe that for a rational function with a pole of order M at $a$ the corresponding M basis functions have the form 
\begin{equation*}
\phi_n(x) = \e^{\mi{n}\theta(x)} \frac{1}{x-\cnj{a}_n}\qquad (n=1..M).
\end{equation*}

These are localized "Fourier like"  basis functions around the real part of  $a$ scaled by the imaginary part.

\subsection*{ Example of a multiscale Wavelet decomposition}

The infinite  Blaschke products
\begin{equation*}
G_n(x) = \prod_{j\le n} \frac{j-\mi}{j+\mi}\, \frac{x-j-\mi}{x-j+\mi}
\text{\quad and\quad } G(x) = \prod_{j\in {\mathbb Z}} \frac{j-\mi}{j+\mi}\,
\frac{x-j-\mi}{x-j+\mi}
\end{equation*}
 can be expressed in terms of known functions:
\begin{equation}\label{G}
G_n(x) = \frac{\Gamma(-\mi-n)}{\Gamma(\mi-n)}\,
\frac{\Gamma(x-n+\mi)}{\Gamma(x-n-\mi)} \text{\quad and\quad } G(x) =
\frac{\sin \pi(\mi-x)}{\sin \pi(\mi+x)}.
\end{equation}

\subsection*{An orthonormal system}

Consider the function
$\phi(x) = \displaystyle\frac{\Gamma(x-1+\mi)}{\sqrt{\pi}\Gamma(x-\mi)}$.
It is easily checked that
$$\phi(x-n) =
\frac{\Gamma(\mi-n)}{\Gamma(-\mi-n)} \, \frac{G_n(x)}{\sqrt{\pi}\bigl(
  x-(n+1)+\mi\bigr)}.$$
Set $\phi_n(x) = \phi(x-n)$. For fixed~$m$, the functions
$\phi_n/G_m$, for $n\ge m$, form a Malmquist-Takenaka basis of
$(G/G_m)\hardy{2}$. In other terms, the functions $\phi_n$, for $n\ge
m$, form an orthonormal basis of $G_m\hardy{2}\ominus
G\hardy{2}$. This means that the functions $\phi_n$ (for $n\in
{\mathbb Z}$) form a Malmquist-Takenaka basis of the orthogonal
complement of $G\hardy{2}$ in $\hardy{2}$.

\subsubsection*{Multiscale decomposition}

As $|1-G(2^nx)|\le C2^{n}$ all the products
$\displaystyle {\mathscr B}_n(x) = \prod_{j< n} G(2^jx)$
are convergent and $\displaystyle \lim_{n\to -\infty} {\mathscr B}_n =
1$ uniformly.

Let ${\mathscr B}={\mathscr B}_0$. Obviously, ${\mathscr B}_n(x) =
{\mathscr B}(2^nx)$.
Consider the following subspaces of $\hardy{2}$:
${\mathsf E}_n = {\mathscr B}_n\hardy{2}$.
This is a decreasing sequence. The space $\displaystyle {\mathsf
  E}_{+\infty} = \bigcap_{n\in {\mathbb Z}} {\mathsf E}_n$ is equal
to~$\{0\}$ since a function orthogonal to this space would have too
many zeros, and the space $\displaystyle {\mathsf E}_{-\infty} =
\mathrm{closure~of}\bigcup_{n\in {\mathbb Z}} {\mathsf E}_n$ is
equal to $\hardy{2}$ since ${\mathscr B}_n$ converges uniformly to 1
when $n$ goes to $-\infty$.

For all $n$ and~$j$, let
$$\phi_{n,j}(x) =
2^{n/2}\phi(2^nx-j){\mathscr B}(2^{n}x).$$
Then, for all~$n$, $(\phi_{n,j})_{j\in {\mathbb Z}}$ is an orthonormal
basis of ${\mathsf E}_{n}\ominus {\mathsf E}_{n+1}$. We conclude that
$(\phi_{n,j})_{n,j\in{\mathbb Z}}$ is an orthonormal basis of
$\hardy{2}$.

\section{ Adapted MT bases, "phase unwinding"}
We now  find  a "best"  adapted Malmquist Takenaka basis to analyze or unwind the oscillations of  a given function.

The idea is to peel off the oscillation of a function by dividing by its  Blaschke product defined by the zeroes of the function, this procedure is iterated to yield an expansion in an orthogonal collection of  functions or Blaschke products which of course are naturally embedded in a MT basis, once the zeroes are ordered.

\subsection*{The unwinding series.} There is a natural way to  iterate the Blaschke factorization, it is inspired by the power series expansion
of a holomorphic function on the disk. If $G$ has no zeroes inside $\mathbb{D}$, its Blaschke factorization is the trivial one $G = 1 \cdot G$, however,
the function $G(z)-G(0)$ certainly has at least one root inside the unit disk $\mathbb{D}$ and will therefore yield some nontrivial Blaschke factorization
$G(z) - G(0) = \blp_1 G_1$.
We write
\begin{eqnarray*}
 F(z) &=& \blp(z) \cdot G(z)
 = \blp(z) \cdot \bigl(G(0) + (G(z) - G(0)\bigr) \\
&=& \blp(z) \cdot \bigl(G(0) + \blp_1(z) G_1(z)\bigr) 
= G(0) \blp z + \blp(z) \blp_1(z) G_1(z).
\end{eqnarray*}
An iterative application gives rise to  the \textit{unwinding series}
$$ F = a_1 \blp_1 + a_2 \blp_1 \blp_2 + a_3 \blp_1 \blp_2 \blp_3 + a_4 \blp_1 \blp_2 \blp_3 \blp_4+ \dots$$
This orthogonal expansion first appeared in the PhD thesis of Michel Nahon~\cite{nahon} and independently by T.~Qian in \cite{qtao,qw} Detailed approximations in smoothness spaces were derived by S. Steinerberger in \cite {coifman}. Given a general function $F$ it is not numerically feasible to actually compute the roots of the function; a crucial insight in \cite{nahon} is that this is not necessary -- one can numerically obtain the Blaschke product
in a stable way by using a method that was first mentioned in a paper of Guido and Mary Weiss \cite{ww} and has been investigated with respect to stability by Nahon \cite{nahon} Using the boundedness of the Hilbert transform one can prove easily convergence in $L^p, 1<p<\infty$.
 
\subsection* { The fast algorithm of  Guido and Mary Weiss \cite{ww}  }

Our starting point is the theorem that any Hardy function   can be
decomposed as
$$ F = \blp \cdot G,$$
where $B$ is a Blaschke product, that is a function of the form
$$ \blp(z) = z^m\prod_{i \in I}{\frac{\overline{a_i}}{|a_i|}\frac{z-a_i}{1-\overline{a_i}z}},$$
where $m \in \mathbb{N}_{0}$ and $a_1, a_2, \dots \in \mathbb{D}$ are zeroes inside the unit disk $\mathbb{D}$ 
and $G$ has no roots in $\mathbb{D}$. For $|z|=1$ we have $|\blp(z)| = 1$
which motivates the analogy 
$$\blp \sim \mbox{frequency and}~G \sim \mbox{amplitude}$$
for the function restricted to the boundary. However, the function $G$ need not be constant: it can be any function that never vanishes
inside the unit disk. If $F$ has roots inside the unit disk, then the Blaschke factorization $F = \blp\cdot G$ is going to be nontrivial (meaning
$\blp \not\equiv 1$ and $G \not\equiv F$). $G$ should be 'simpler' than $F$ because the winding number around the origin decreases.

In fact since $|F|=|G|$ and $\ln(G)$ is analytic in the disk we have formally that
$G=\exp(\ln|F|+ \mi(\ln|F|)^\sim))=\exp(\mathscr H(\ln|F|))$   where $\mathscr H$ is the projection onto the Hardy space.
and   $\blp=F/G$. 
G can be computed easily using the FFT \cite{nahon}.

\subsection*{ A remarkable unwinding}
The following is an explicit unwinding of a singular inner function in the upper half plane  illustrating this exponentially fast  approximation
of  $\exp \frac{2\mi\pi}{x}$:
 \begin{equation*}
\exp \frac{2\mi\pi}{x} = \e^{-2\pi} + \bigl( 1-\e^{-4\pi}\bigr) \sum_{n\ge 0}
(-1)^n\e^{-2n\pi} B(x)^{n+1},
 \end{equation*}
 where $B$ is a Blaschke product whose zeros are $\{1/(j+\mi)\}_{j\in
   {\mathbb Z}}$.

\section{Geometric function theory:  the role of compositions of Blaschke products.}
 \subsubsection*{Iteration of Blaschke products}
 
 We claim that by building Blaschke product through composition  we open up rich dynamical structures, and libraries of corresponding Malmquist Takenaka bases.  

 We are interested in iteration of finite Blaschke products
 \begin{equation*}
\blp(z) = \e^{\mi\theta}z^\nu \prod_{j=1}^\mu \frac{z+a_j}{1+\cnj{a}_jz}, 
\end{equation*}
where~$\mu$ and~$\nu$ are nonnegative integers and the $a_j$ are complex
numbers of modulus less than~1.

It is well known that ${\mathbb T}$ and ${\mathbb D}$ are globally
invariant under~$\blp$, as well as the complement of~$\cnj{\mathbb D}$ in
the Riemann sphere.

A careful discussion can be found in~\cite{CP}. Here is the main result.
\begin{theorem}\label{Diverge}
Let $\blp$ be a finite Blaschke product with a fixed point~$\alpha$
inside the unit disk. Then there exists a sequence
$\alpha,a_1,a_2,\dots,a_j,\dots$ of complex numbers in the unit disk
and an increasing sequence $(\nu_j)_{j\ge1}$ of positive integers such
that $a_1,a_2,\dots,a_{\nu_n}$ are the zeros, counted according to
their multiplicity, of~$\blp_{n}$ (the nth iterate of $\blp$). Moreover $\displaystyle \sum_{j\ge  1}(1-|a_j|) = +\infty$. Also, $\blp_n$ converges towards~$\alpha$ unformly on compact subsets of the open unit disk.
\end{theorem}

\subsection*{ Dynamic Multiscale analysis through composition of Blaschke products }

Each Blaschke product $\blp$ defines invariant subspaces of
$\mat{H}^p$. The projection on this space is given by the kernel
$\displaystyle \frac{\blp(z)\cnj{\blp(w)}}{z-w}$. This projection is
continuous for $1< p< +\infty$.

Let $F$ be a Blaschke product of degree at least~2 with a fixed point
inside the unit disk. Its iterates define a hierarchy of nested
invariant subspaces ${\mathsf E}_n = F_n\mat{H}^2$.

Due to Theorem~\ref{Diverge},
$\displaystyle\bigcap_{n\ge 1} {\mathsf E}_n = \{0\}$.

The Takenaka construction provides orthonormal bases of $E_n\ominus
E_{n+1}$. But this is not canonical as it depends on an ordering of the
zeros of $F_{n+1}/F_n$.

Figure~1 shows 1st, 3rd, and 5th iterates of $F(z)
=z(z-2^{-1})/(1-2^{-1}z)$. Figure~2 displays the phase for the fourth  iterate of $F(z)
=z^{2}(z-2^{-2})/(1-2^{-2}z)$. The upper pictures display the phases
modulo $2\pi$ (values in the interval $(-\pi,\pi]$) of theses Blaschke
products while the lower pictures display minus the logarithms of
their absolute value. The coordinates $(x,y)$ correspond to the point
$\e^{-y + \mi x}$. On these figures it is easy to locate the zeros,
specially by looking at the phase which then has an abrupt jump.

\subsection*{Remarks on Iteration of Blaschke products as a "Deep Neural Net"}
 In the upper half plane let
$(a_j)_{1\le j}$ be a finite sequence  of complex numbers
with positive imaginary parts.
The corresponding Blaschke product on the line is
$$\blp(x) = \prod_{j\ge 0} \frac{x-a_j}{x-\cnj{a}_j}.$$

 We can write $\blp(x)=\exp\bigl(\mi\theta(x)\bigr)$, where 
$$\theta(x)=\sum_{j\ge 0}\sigma\bigl((x-\alpha_j)/\beta_j\bigr) $$
with $a_j=\alpha_j+\mi\beta_j$  and $\sigma=\arctan x + \pi/2 $ is a sigmoid.

This is precisely the form of a single layer in a Neural Net, each unit has a weight and bias determined by $a_j$. We obtain the various layers of a deep net through the composition 
of each layer with a preceding layer. In our preceding examples we took a single short layer given  by a simple  Blaschke term with two zeroes in the first layer that we iterated to obtain an orthonormal Malmquist Takenaka basis ( we could have composed  different elementary  products at each layer), demonstrating the versatility of the method to generate highly complex functional representations.

As an example let $F(z)$ be mapped from G, \eqref{G} in the section on wavelet construction.
$$ F(z)=G(w)= \frac{\sin(\pi(\mi-w))}{\sin(\pi(\mi+w))}
\text{\quad with\quad} w=\frac{\mi(1-z)}{(1+z)}.$$  

We can view the phase of F as a neural layer which when composed with itself results in a phase which is a two layer neural net represented graphically in fig~3. 

Where each end of a color droplet corresponds to one zero or unit of the two layer net.

We refer to Daubechies et al.\ \cite{ReluDNN} for a description of a similar iteration for piecewise affine functions in which simple affine functions play the role of a Blaschke product. 

\section{ Higher dimensions,  $\theta$-holomorphy}
Our goal is to explore methodologies to use the remarkable analytic approximation theorems described above to enable deeper understanding of real  analysis, in higher dimensions. We  know that Blaschke factorization do not exist, nevertheless there are remarkable bases that can be lifted.

We start by observing that  $Z_{\theta}=(x\cdot\theta+\mi y)=t+\mi y$ is harmonic in $(x,y)$,(in 3 dimensions) and so is $Z_{\theta}^k$
This is a harmonic homogeneous polynomial of degree k in $(x_1,x_2,y)$ that is constant in the direction perpendicular to $\theta$.
here we identified $\theta$ with   $(\cos\theta,\sin\theta) $.
It is well known  \cite{CW} that 

\begin{equation*}
 \frac{1}{2\pi} \int_0^{2\pi}
\e^{-\mi{l}\theta}{Z_{\theta}^k} \,\dif \theta= Y_l^{k}(x_1,x_2,y) \quad 
 (\abs{l}\leq{k})
\end{equation*}
is the standard orthogonal basis of spherical Harmonics in 3 dimensions.

As a consequence we see that any Harmonic function $U(x,y)$ is a superposition over $\theta $ of holomorphic functions in $Z_{\theta}$,
more explicitly a Power series in $Z_{\theta}$ with coefficients depending on  $\theta $.

\begin{equation*}
 U(x,y) = \frac{1}{2\pi} \int_0^{2\pi}
\ F_{\theta} (Z_{\theta)}) \,\dif\theta = \frac{1}{2\pi} \int_0^{2\pi}
{\sum_{k\ge 0} a_{k}({\theta}){Z_{\theta}^k\dif\theta }}
\ 
\end{equation*}  where $a_{k}({\theta})$ is a trigonometric polynomial of degree $k$.

Another   example, taking  $$ F_{\theta} (Z_{\theta})= \e^{-2\mi\pi rZ_{\theta}}{\e^{-2\mi\pi k{\theta}}}.$$  

we get the harmonic function   $$ J_{k}\left(2\pi r\sqrt{x_{1}^2+x_{2}^2}\right) \e^{-2\mi\pi k{\phi} }\e^{-2\pi yr} $$

\section*{ Radon and Fourier in the upper Half space}

This relationship between holomorphic functions in planes  as generating all harmonic functions can most easily be explored through Fourier analysis.
We define the Radon transform, and relate it to the Fourier transform to lead to the $\theta$-holomorphic version.

Let
\begin{equation}
  \radon\theta f(t) = \int_{x\in \theta^\perp} f(t\theta+x)\,\dif x.
\end{equation}
Obviously $\radon{-\theta}f(t) = \radon\theta f(-t).$ The formula $\widehat{\radon\theta f}(t) = \hat{f}(u\theta)$ for Fourier transforms is well known.

For $f\in L^1({\mathbb R}^n)$, consider its harmonic extension~$u$ to ${\mathbb R}_+^{n+1}$. For $x\in {\mathbb R}^n$ and $y>0$ we have
\begin{eqnarray*}
  u(x,y) &=& f\star {\mathsf P}_y(x) = \int\e^{2\mi\pi x\cdot\xi}\e^{-2\pi|\xi|y}\hat{f}(\xi)\,\dif \xi\\
         &=& \int_{S_{n-1}}\left(\int_0^\infty \e^{2\mi\pi r(x\cdot\theta+\mi y)} \hat{f}(r\theta)r^{n-1}\dif r\right)\,\dif\theta\\
         &=& \int_{S_{n-1}} F_\theta(x\cdot\theta+\mi y)\,\dif\theta,
\end{eqnarray*}
where
\begin{equation}\label{F}
  F_\theta(z) = \int_0^\infty \e^{2\mi\pi r z} \hat{f}(r\theta)r^{n-1} \dif r
 = \int_0^\infty \e^{2\mi\pi r z} \widehat{\radon\theta f}(r)\,r^{n-1}\dif r.
\end{equation}

When $n=2$, we have 
\begin{equation*}
  \widehat{F_\theta}(t) = \widehat{\radon\theta f}(t)\, t \ind_{(0,+\infty)}(t)
                        = \frac{1}{2\mi\pi}\widehat{\Dif\radon\theta f}(t)\,\ind_{(0,+\infty)}(t).
\end{equation*}
So, for $\Im z>0$,
\begin{equation*}
  F_\theta(z)  = -\frac{1}{4\pi^2} \int_{-\infty}^{\infty} \frac{\dif(\radon\theta f(t))}{t-z}\\
        =    -\frac{1}{4\pi^2} \int_{-\infty}^{+\infty} \frac{\radon\theta f(t)}{(t-z)^2}\,\dif t.
\end{equation*}

For general~$n$ we get
$$ F_\theta(z) = \frac{(n-1)!}{(2\mi\pi)^{n}} \int_{-\infty}^{+\infty} \frac{\radon\theta f(t)}{(t-z)^{n}}\,\dif t.
$$

\subsection*{Some isometries}

We describe some computations in the case when $n=2$ and mention the case of higher dimension at the end of this section.

In view of~\eqref{F}
\begin{equation}\label{FT}
 \widehat{F_\theta(\cdot+\mi y)}(r) = \hat{f}(r\theta)\,\e^{-2\pi r y}r\ind_{(0,\infty)}(r).
\end{equation}
Hence, the Plancherel identity yields
\begin{eqnarray*}
  \int_0^\infty \dif y \int_{-\infty}^{+\infty}|F_\theta(t+\mi y)|^2\dif t &=& \int_0^\infty \int_0^\infty |\hat{f}(r\theta)|^2 r^2\e^{-4\pi r y}\dif r \dif y\\
&=& \frac{1}{4\pi} \int_0^\infty |\hat{f}(r\theta)|^2 r\dif r
\end{eqnarray*}
\bigskip

Let    $\displaystyle \|F\|_{B}^2 = \int_0^\infty \int_{-\infty}^{+\infty}|F_\theta(t+\mi y)|^2\dif y\dif t$ (this is the norm of the Bergman space on the upper half plane). Then
\begin{equation}\label{bergman}
  4\pi\int_0^{2\pi}\|F_\theta\|_{B}^2\dif\theta = \iint_{(0,+\infty)\times (0.2\pi)} \left|\hat{f}(r\theta)\right|^2 r\,\dif r \dif\theta
            = \|f\|_{L^2({\mathbb R}^n)}^2.
\end{equation}
\medskip

\noindent We have 
\begin{equation*}
  \frac{\partial u(x,y)}{\partial y} = -2\pi \int_{{\mathbb R}^2} \e^{2\mi\pi \xi\cdot x} |\xi| \e^{-2\pi|\xi|y}\hat{f}(\xi)\,\dif \xi.
\end{equation*}

\begin{eqnarray}
  \iint_{{\mathbb R}_+^3} \left|\frac{\partial u(x,y)}{\partial y}\right|^2\dif x\dif y
     &=& (2\pi)^{2} \int_{{\mathbb R}^2} \left(\int_0^\infty \e^{-4\pi|\xi| y} \dif y\right)|\hat{f}(\xi)|^2|\xi|^{2}\dif \xi\nonumber\\
&=& \pi\int_{{\mathbb R}^n}|\hat{f}(\xi)|^2 |\xi|\dif \xi\label{deriv}
\end{eqnarray}
Equation~\eqref{FT} yields
  $\displaystyle \int_{-\infty}^{+\infty} |F_\theta(t)|^2\,\dif t = \int_0^\infty |\hat{f}(r\theta)|^2 r^2\dif r$,
and
\begin{equation}\label{H2}
  \int_0^{2\pi} \|F_\theta\|_{\hardy2({\mathbb R})}^2 \dif \theta = \int_{{\mathbb R}^2} |\hat{f}(\xi)|^2 |\xi|\dif \xi.
\end{equation}
Formulas~\eqref{deriv} and~\eqref{H2} together give
\begin{equation}\label{dirichlet}
  \int_0^{2\pi} \|F_\theta\|_{\hardy2({\mathbb R})}^2 \dif \theta = \frac1\pi \iint_{{\mathbb R}_+^{3}} \left|\frac{\partial u(x,y)}{\partial y}\right|^2\dif x\dif y,
\end{equation}

\medskip

In higher dimension formulas~\eqref{bergman} and~\eqref{dirichlet} become
$$
\int_{S_{n-1}} \dif \theta \iint_{{\mathbb R}_+^{n+1}} |F_\theta(t+\mi y)|^2y^{n-2}\dif t\dif y = \frac{1}{(4\pi)^{n-1}} |f\|_{L^2({\mathbb R}^n)}^2$$
and
\begin{multline}\label{dirichlet2}
\iint_{{\mathbb R}_+^{n+1}} \left|\frac{\partial^k u(x,y)}{\partial y^k}
\right |y^{2k-n}\dif x\dif y\\ = \frac{(4\pi)^{n-1}\Gamma(2k-n+1)}{2^{2k}} \int_{S_{n-1}} \|F_\theta\|_{H^2(\mathbb R)}^2\dif\theta.
\end{multline}

\subsection*{Remarks;   "lifted Analysis" of Harmonic functions}
One of our goals is to enable the application of some of the one dimensional analytic approximation tools to higher dimensions. We refer to Michel Nahon's thesis \cite{nahon} in which he decomposes a function in the plane as a sum of functions whose Fourier transform live in  thin wedges, as a tool to extract features ( gradients of Phase) from an image of a fingerprint.
This illustrates potential variants of our current approach.

We envision a function in two variables  represented as a superposition of $F_{\theta} (t + iy)$,  each of which is approximated to  error $\epsilon$ leading to a harmonic function approximation of error  $\epsilon$   in the Dirichlet space. Similar estimations with different mix of Hilbert spaces can be easily derived as in~\cite{coifman}, leading to faster rates of convergence (when more regularity is present).

Another obvious application is the representation of a Calder\'on Zygmund operator given by  a Fourier multiplier homogeneous of degree $0$,   $\Omega(\theta)$, simply by averaging 
 $ \Omega(\theta)  F_{\theta} (t + iy)$.
 
 The representation of these operators is a version of the rotation  method (no parity required on $\Omega$ ). Also it provides a local representation method for  generalized conjugate functions  or C-Z operators, just by using the local spherical Harmonics version of the ${\theta}$-holomorphic representation. 
 In particular Riesz transforms correspond to   $\Omega(\theta)=(\cos\theta,\sin\theta) $. 
 
 \subsection*{A natural ortho-basis in the Dirichlet space}
 
 We now use identity~\eqref{dirichlet} to transfer an orthonormal basis of the Hardy space~$H^2$ to an orthonormal system in the Dirichlet space in ${\mathbb R}_+^{3}$.
 
 Start from the basis $\frac{\mi}{\sqrt{\pi}} \left(\frac{z-i}{z+i}\right)^n\frac{1}{z+i}$ of $H^2$ ( this corresponds to the Fourier basis in the  disc, mapped  to the upper half plane Hardy space). We consider the generating function
$$F(z) = \frac{\mi}{\sqrt{2\pi^3}}\sum_{k\ge 0} (z-i)^nt^n/(z+i)^{n+1}$$ and compute  $G(x,y) = \displaystyle \int_0^{2\pi} F(z_\theta)\, \dif \theta$.
%
%
We get
$$G(x,y) = \sqrt{2/\pi}/\sqrt{(\rho^2+2y+1)t^2-2(1-\rho^2)t+\rho^2-2y+1},$$
where $\rho = \sqrt{x_1^2+x_2^2+(y+1)^2}$.\bigskip

This also can be written as

$$G(x,y) = \frac{\sqrt{2/\pi}}{\sqrt{\rho^{2}+2 y +1}\, \sqrt{a^{2} t^{2}-2 b a
    t +1}},$$  if one sets
$a = \sqrt{\frac{\rho^{2}-2 y +1}{\rho^{2}+2 y +1}}$ and
$b = \frac{\rho^2-1}{\sqrt{(\rho^2+1)^2-4y^2}}.$
\medskip

It results that the functions
$\displaystyle \frac{\sqrt{2/\pi}a^n P_{n}(b)}{\sqrt{\rho^{2}+2
    y +1}}$,
where the $P_n$ are the Legendre polynomials, form an orthonormal system in the Dirichlet space.
\medskip

To get an orthonormal basis of the Dirichlet space in 3 dimensions, it suffices to take
$\displaystyle \frac{\sqrt{2/\pi}2a^n P_{n}(b)\,{\mathrm e}^{\mi k \theta}}{\sqrt{\rho^{2}+2
    y +1}}$,
with $k\in \mathbb Z$ and $n\ge 0$.

Of course such computation can be done in higher dimension: isometry~\eqref{dirichlet2} allows to transfer orthonormal bases of $H^2$ to orthonormal systems in a suitable Dirichlet space.

\subsection*{ Concluding remarks and potential applications}

As we all know complex methods, such as interpolation of operators, or the remarkable proofs by Calder\'on of the boundedness in $L^2$ of commutators with the Hilbert transform, or the Cauchy integral on Lipschitz curves are powerful tools. Over the years the goal has been to convert them into real variable methods. In parallel the quest for higher dimensional complex tools is continuing, see the examples  \cite{CW} in which various systems generalizing holomorphic functions to higher dimension are studied. The point here, is that the infinite dimensional $\theta$-holomorphic functions generate all of these systems through the choice of appropriate multipliers (as described for the Riesz system) .

Our goal here was to describe recent nonlinear  analytic tools in the classical setting. and transfer  them to the higher dimensional real setting.
Together with  Guido Weiss we had observed~\cite{CW} that all harmonic functions in higher dimensions are combinations of holomorphic functions on subplanes which are constant in normal directions. The recent developments in one dimension as well as the isometries  described here, and the corresponding efficient approximation methods, open the door for applications in higher dimensions, such as image denoising. See \cite{CSW}  for the impact of  unwinding on  precision Doppler analysis in 1 dimension, which we expect to carry over to 2 or 3 dimensions.

Observe also that, for simplicity,  we restricted our attention to 2 dimensions in cylindrical coordinates. We could have defined more generally power series in  the variable   
$ Z_{\epsilon}= (x\cdot{\epsilon)} $    
where ${\epsilon}$  satisfying;   $(\epsilon\cdot\epsilon)=0$ ,  represents a point on the complex quadric with $|\Re{\epsilon}|=1,|\Im{\epsilon}|=1$. , or a two dimensional plane spanned by $\Re{\epsilon},\Im{\epsilon}$.

Clearly we can extend the preceding discussion to this setting.  Where; $ Z_{\epsilon}= (x\cdot{\epsilon)} $ is the point $ t+\mi s$ in the complex plane with  coordinate    
$t\Re{\epsilon}+\mi s\Im{\epsilon}$.  

\begin{figure}
  \includegraphics[width=115mm,height=60mm]{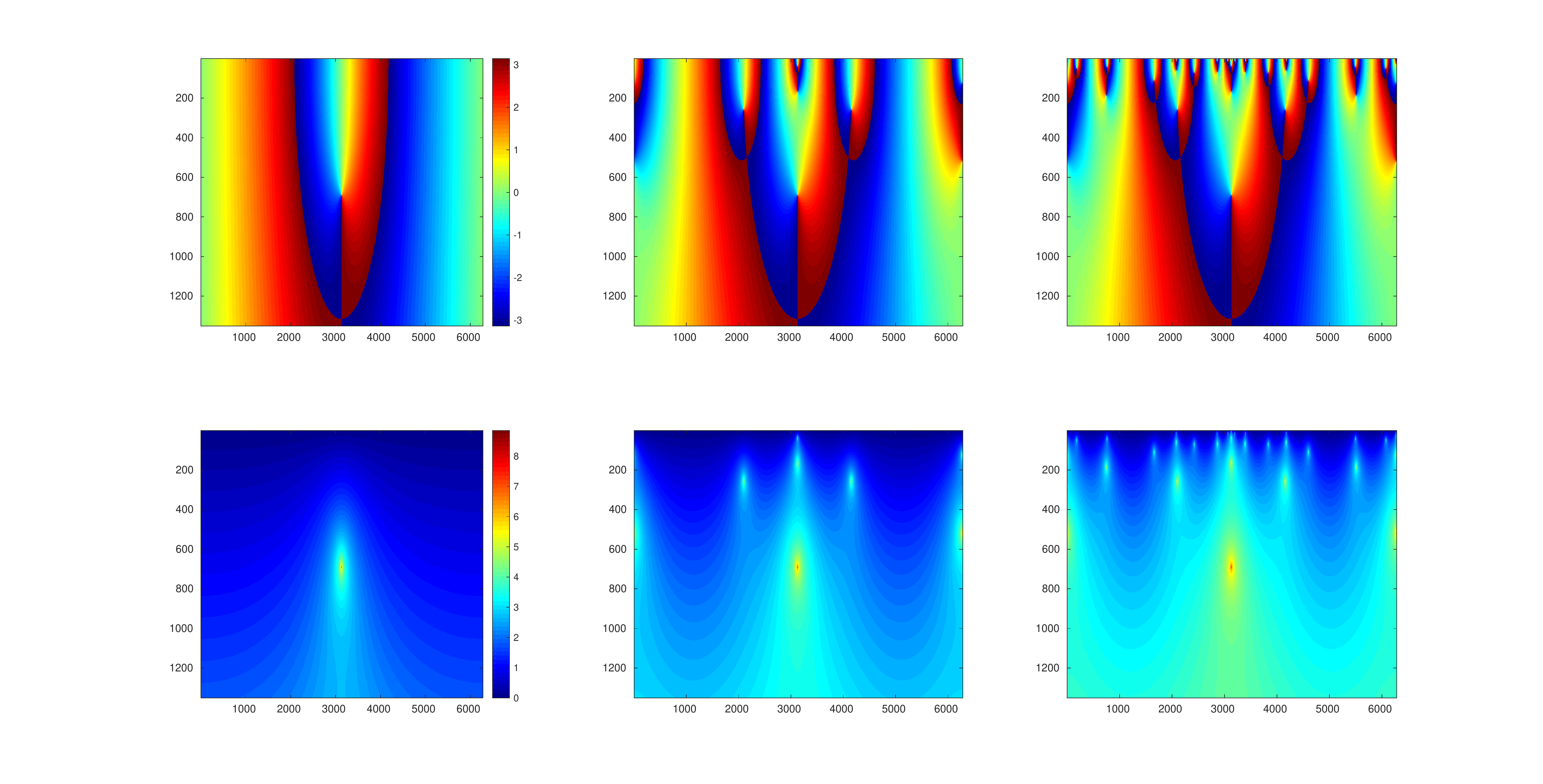}\\
\caption{ The argument and the absolute value of $F(z)$, $F^{(3)}(z)$,
  and $F^{(5}(z)$,\\ with $F(z)=\frac{z(z-2^{-1})}{1-2^{-1}z}$ and $z=
  \exp(-y+\mi x)$.}
\end{figure}

\begin{center}
\begin{figure}
\includegraphics[width=115mm,height=60mm]{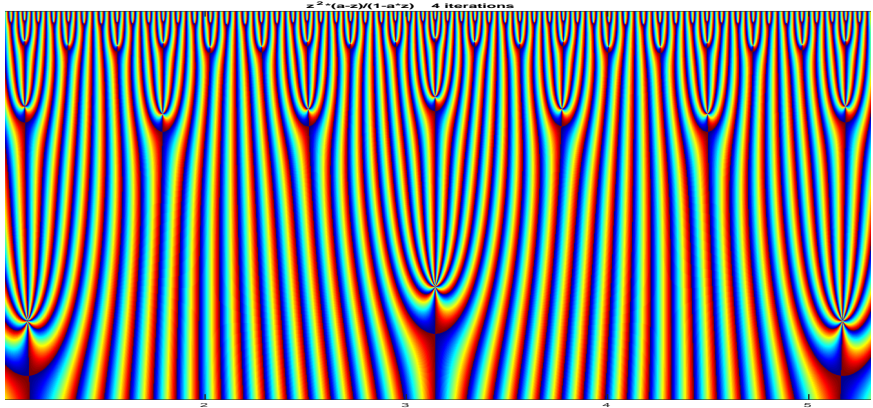}
\caption {The multiscale view of the argument of $F^{(4)}(z)$,
 with $F(z)=\frac{z^2(z-2^{-2})}{1-2^{-2} z}$    
 and $z=\exp(-y+\mi x)$.}
\end{figure}
\end{center}

\begin{figure}
\includegraphics[width=115mm,height=70mm]{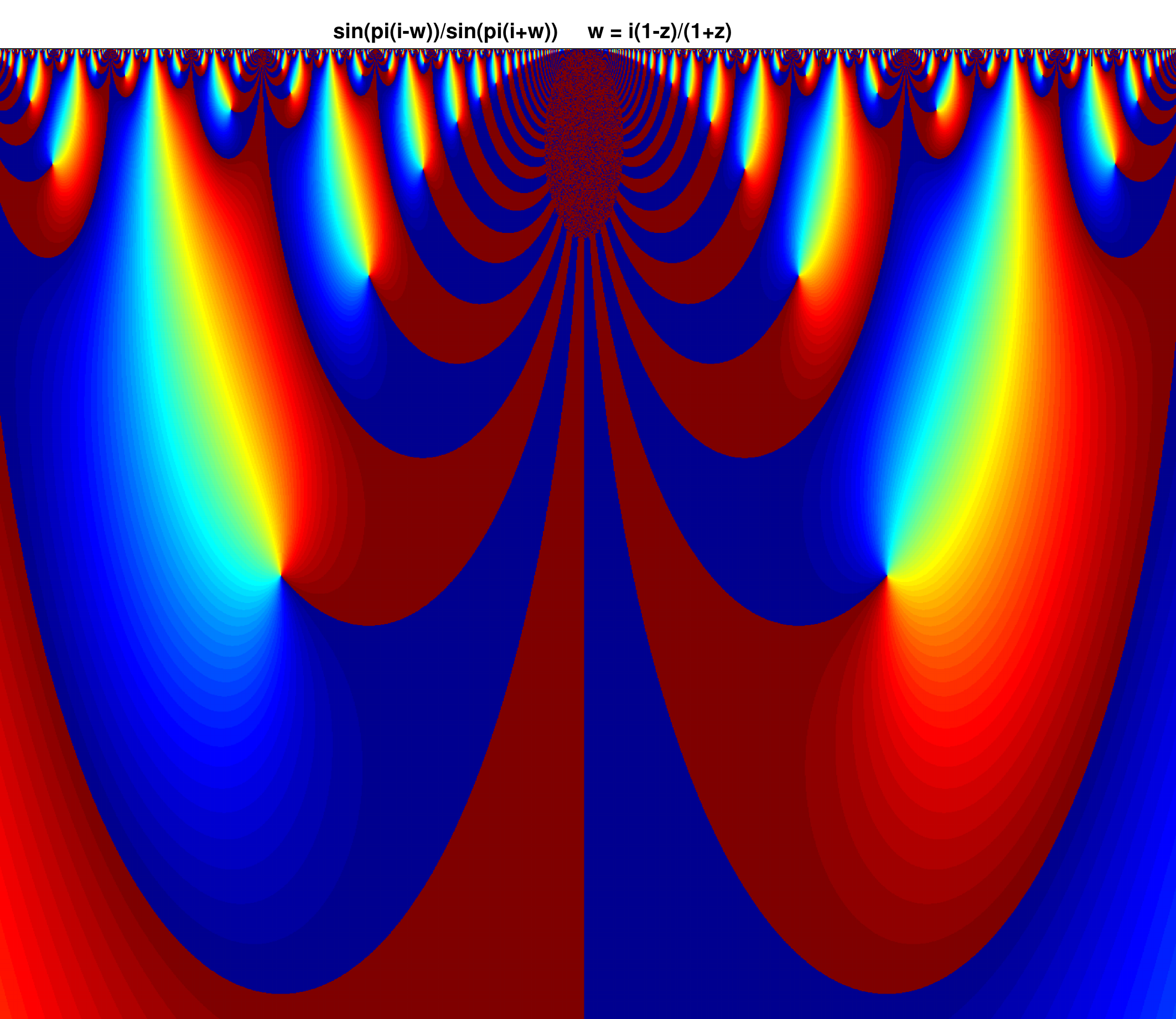}
\caption{ Two iterations  of $ F(z)=G(w)= \frac{\sin(\pi(\mi-w))}{\sin(\pi(\mi+w))}
\text{\quad with\quad} w=\frac{\mi(1-z)}{(1+z)}.$ }  

\end{figure}


\begin{thebibliography}{99}

\bibitem{brolin} Brolin,~H., Invariant sets under iteration of
  rational functions, \textsl{Arkiv f\"or Matematik, \textbf{6-6},
    (1965),~103--144.}
  
  Springer 1993.

\bibitem{CP}  Coifman,~R.~R., and Peyri\`ere,~J., Phase Unwinding, or invariant subspace decompositions of Hardy Spaces. \textsl{Journal of Fourier Analysis and Applications}  \textbf{25} (2019),~684--695.

\bibitem{CW1}  Coifman,~R.~R., and Weiss,~G., A kernel associated with certain multiply connected domains, and its applications to factorization theorems.  \textsl{ Studia Mathematica\ (1966)}.
  
\bibitem{coifman} Coifman,~R.~R., and Steinerberger,~S., Nonlinear
  phase unwinding of functions. \textsl{J.~Fourier
    Anal.\ Appl.\ (2016)},~1--32.

\bibitem{CSW} Coifman,~R.~R., Steinerberger,~S., and Wu,~H.~T., Carrier
  frequencies, holomorphy and unwinding. arXiv preprint
  arXiv:1606.06475, 2016 - arxiv.org.
  
\bibitem{CW} Coifman,~R.~R., Weiss,~G., Analyse Harmonique Non-Commutative sur Certains Espaces Homogenes, Lecture Notes i Mathematics 242, Springer-Verlag.

\bibitem{ReluDNN} I. Daubechies, R. DeVore, S. Foucart, B. Hanin, and G. Petrova.Nonlinear Approximation and (Deep) ReLU Networks.
arXiv:1905.02199v1  [cs.LG]  5 May 2019.


\bibitem{feichtinger} Feichtinger, H.G. and Pap, M., Hyperbolic
  wavelets and multiresolution in the Hardy space of the upper half
  plane, \textsl{Blaschke Products and Their Applications, (2013),}
  Springer.  



  
 
 \bibitem{malmquist} Malmquist,~F,, Sur la determination d’une classe
de fonctions analytiques par leurs valeurs dans un
ensemble donne de poits, \textsl{C.R. 6ieme Cong.\
Math.\ Scand.\ (Kopenhagen, 1925), Copenhagen,} (1926),
Gjellerups,~253--259.
 
 \bibitem{qtao} Mi,W., Qian,~T., and Wan,~F., A Fast Adaptive Model Reduction Method Based on Takenaka-Malmquist Systems, Systems \& Control Letters. Volume 61, Issue 1, January 2012, Pages 223--230.
 
\bibitem{nahon} Nahon, M., Dissertation, Yale University (2000).


  




 \bibitem{qw} Qian,~T., I. T. Ho, Leong~,~I.~T.,  and Wang,~Y.~B., Adaptive decomposition of functions into pieces of non-negative instantaneous frequencies, International Journal of Wavelets, Multiresolution and Information Processing, 8 (2010), no. 5, 813--833.

  

\bibitem{takenaka} Takenaka, S., On the orthogonal functions and a new formula of interpolation,
\textsl{Jpn. J. Math. II } (1925),~129--145.
 
 \bibitem{ww} Weiss,~G,\  and Weiss~M,  A derivation of the main results of the theory of $H^p$-spaces. Rev. Un. Mat. Argentina \textbf{20}  (1962),~63--71.
 \end{thebibliography}
\end{document}